\documentclass{amsart}
\usepackage{eurosym}
\usepackage{amsfonts}
\usepackage{amsmath}
\usepackage{amssymb}
\usepackage{graphicx}
\usepackage{amscd}

\newtheorem{theorem}{Theorem}
\theoremstyle{plain}

\newtheorem{definition}{Definition}

\newtheorem{lemma}{Lemma}

\newtheorem{proposition}{Proposition}
\newtheorem{remark}{Remark}

\numberwithin{equation}{section}

\begin{document}
\title[A general strong law and Applications]{A general strong law of large
numbers and applications to associated sequences and to extreme value theory}
\author{Harouna SANGARE}
\author{Gane Samb LO}
\address[Harouna SANGARE and Gane Samb LO]{ Departement d'Etudes et Recherches en Mathematiques et d'Informatique, Faculté des Sciences Techniques (FST), Université des Sciences Techniques et Technologies de Bamako (USTTB), Bamako, Mali\\
Laboratoire de Statistique Théorique et Appliquée (LSTA), Université Pierre et Marie Curie, Paris, France}
\email[Harouna SANGARE]{harounasangareusttb@gmail.com}
\email[Gane Samb LO]{ganesamblo@ganesamblo.net}

\begin{abstract}
The purpose of this paper is to establish a general strong law of large
numbers (SLLN) for arbitrary sequences of random variables (rv's) based on
the squared indice method and to provide applications to SLLN of associated
sequences. This SLLN is compared to those based on the H\'{a}jek-R\'{e}nyi
type inequality. Nontrivial examples are given. An interesting issue that is
related to extreme value theory (EVT) is handled here.
\end{abstract}

\maketitle



%


\Large
\noindent

\section{Introduction}

\noindent \label{sec0} In this paper, we present a general SLLN for
arbitrary rv's and particularize it for associated sequences. In the recent
decades both strong law of large numbers and central limit theorem for
associated sequences have received and are still receiving huge interests
since Lebowitz \cite{lebow} and Newman \cite{newman84} results under the
strict stationarity assumption. The stationarity assumption was dropped by
Birkel \cite{birkel}, who proved a version of a SLLN that can be interpreted
as a generalized Kolmogorov's one. A recent account of such researches in
this topic is available in \cite{rao}. Although many results are available
for such sequences, there are still many open problems, especially regarding
nonstationary sequences.

\bigskip

\noindent We intend to provide a more general SLLN for associated sequences
as applications of a new general SLLN for arbitrary rv's. This new general
SLLN is used to solve a remarkable issue of extreme value theory by using a
pure probabilistic method.

\bigskip

\noindent Here is how this paper is organized. Since association is the
central notion used here, we first make a quick reminder of it in section %
\ref{subsec1}. In section \ref{subsec2}, we make a round up of SLLN's
available in the literature with the aim of comparing them to our findings.
In section \ref{sec1}, we state our general SLLN for arbitrary rv's and
derive some classical cases. In section \ref{sec2}, we give an application
to EVT where the continuous Hill's estimator is studied by our method. The
section \ref{conclu} concerns the conclusion and some perspectives are
given. The paper is ended by the Appendix section, where are postponed the
proofs of Propositions \ref{prop2} and \ref{prop3} stated in section \ref%
{sec2}.

\bigskip

\noindent To begin with, we give a short reminder of the concept of
association.

\section{A brief reminder of the concept of association}

\label{subsec1} \noindent The notion of positive dependence for random
variables was introduced by Lehmann (1966) (see \cite{lehmann}) in the
bivariate case. Later this idea was extended to multivariate distributions
by Esary, Proschan and Walkup (1967) (see \cite{esary}) under the name of
association. The concept of association for rv's generalizes that of
independence and seems to model a great variety of stochastic models. This
property also arises in Physics, and is quoted under the name of FKG
property (Fortuin, Kastelyn and Ginibre (1971), see \cite{fortuin}), in
percolation theory and even in Finance (see \cite{jiazhu}). \noindent The
definite definition is given by Esary, Proschan and Walkup (1967) (see \cite%
{esary}) as follows.

\begin{definition}
A finite sequence of random variables $(X_{1},...,X_{n})$ is associated if
for any couple of real and coordinate-wise non-decreasing functions $f$ and $%
g$ defined on $\mathbb{R}^{n}$, we have 
\begin{equation}
Cov(f(X_{1},...,X_{n}),\ g(X_{1},...,X_{n}))\geq 0  \label{asso}
\end{equation}%
whenever the covariance exists. An infinite sequence of random variables is
associated whenever all its finite subsequences are associated.
\end{definition}

\noindent We have a few number of interesting properties to be found in (%
\cite{rao}) :

\bigskip

\noindent \textbf{(P1)} A sequence of independent rv's is associated.
\noindent \textbf{(P2)} Partial sums of associated rv's are associated.
\noindent \textbf{(P3)} Order statistics of independent rv's are associated.
\noindent \textbf{(P4)} Non-decreasing functions and non-increasing
functions of associated variables are associated. \noindent \textbf{(P5)}
Let the sequence $Z_{1},Z_{2},...,Z_{n}$ be associated and let $%
(a_{i})_{1\leq i\leq n}$ be positive numbers and $(b_{i})_{1\leq i\leq n}$
real numbers. Then the rv's $a_{i}(Z_{i}-b_{i})$ are associated.\newline

\noindent As immediate other examples of associated sequences, we may cite
Gaussian random vectors with nonnegatively correlated components (see \cite%
{pitt}) and homogeneous Markov chains (see \cite{daley}).\newline

\noindent The negative association was introduced by Joag-Dev and Proschan
(1983) (see \cite{joag}) as follows

\begin{definition}
The variables $X_{1},...,X_{n}$ are negatively associated if, for every pair
of disjoint subsets nonempty $A,$ $B$ of $\{1,...,n\},$ $A=\{i_{1},...,i_{m}%
\}$, $B=\{i_{m+1},...,i_{n}\}$ and for every pair of coordinatewise
nondecreasing functions $f:\mathbb{R}^{m}\rightarrow \mathbb{R}$ and $g:%
\mathbb{R}^{n-m}\rightarrow \mathbb{R}$,%
\begin{equation}
Cov(f(X_{i},\text{ }i\in A),\ g(X_{i},\text{ }i\in B))\leq 0  \label{negasso}
\end{equation}%
whenever the covariance exists. An infinite collection is said to be
negatively associated if every finite sub-collection is negatively
associated.

\begin{remark}
For negatively associated sequences, we have (\ref{negasso}), so the
covariances are non-positive. This remark will be used in sub-subsection \ref%
{subsec6}.
\end{remark}
\end{definition}

\noindent A usefull result of Newman (see \cite{newman80}) on assocation,
that is used in this paper, is the following

\begin{lemma}[Newman \protect\cite{newman80}]
\label{lem1} Suppose that $X$ and $Y$ are two random variables with finite
variance and, $f$ and $g$ are $%
\mathbb{C}
^{1}$ complex valued functions on $\mathbb{R}^{1}$ with bounded derivatives $%
f^{\prime }$ and $g^{\prime }.$ Then 
\begin{equation*}
|Cov(f(X),\text{ }g(Y))|\leq ||f^{\prime }||_{\infty }||g^{\prime
}||_{\infty }Cov(X,\text{ }Y).
\end{equation*}
\end{lemma}

\noindent Here, we point out that strong laws of large numbers and, central
limit theorem and invariance principle for associated rv's are available.
Many of these results in that field are reviewed in \cite{rao}. Such studies
go back to Lebowitz (1972) (see \cite{lebow})\ and Newman (1984) (see \cite%
{newman84}). As Glivenko-classes for the empirical process for associated
data, we may cite Yu (1993) (see \cite{yu93}). We remind the results of such
authors in this :

\begin{theorem}[Lebowitz \protect\cite{lebow} and Newman \protect\cite%
{newman84}]
\label{theo0} Let $X_{1},X_{2},\cdots $be a strictly stationary sequence
which is either associated or negatively associated, and let $T$ denote the
usual shift transformation, defined so that $%
T(f(X_{j_{1}},...,X_{j_{m}}))=f(X_{j_{1}+1},...,X_{j_{m}+1})$. Then $T$ is
ergodic (i.e., every $T$-invariant event in the $\sigma $-field generated by
the $X_{j}$'s has probability $0$ or $1$) if and only if 
\begin{equation}
\lim_{n\rightarrow +\infty }\frac{1}{n}\sum_{j=1}^{n}cov(X_{1},\text{ }%
X_{j})=0.  \label{lebonew}
\end{equation}%
In particular, if \ (\ref{lebonew}) is valid, then for any $f$ such that $%
f(X_{1})$ is $L_{1}$, 
\begin{equation*}
\lim_{n\rightarrow +\infty }\frac{1}{n}\sum_{i=1}^{n}f(X_{i})=\mathbb{E}%
\left( f(X_{1})\right) \ \text{almost surely (}a.s).
\end{equation*}
\end{theorem}

\bigskip \noindent Now we are going to state some classical SLLN's for
arbitrary rv's in relation with H\'{a}jek-R\'{e}nyi's scheme.

\section{Strong laws of large numbers}

\noindent \label{subsec2} For independent rv's, two approaches are mainly
used to get SLLN's. A direct method using squared indice method seems to be
the oldest one. Another one concerns the Kolmogorov's law based on the
maximal inequality of the same name. Many SLLN's for dependent data are
kinds of generalization of these two methods. Particularly, the second
approach that has been developed to become the H\'{a}jek-R\'{e}nyi's method
(see \cite{hajeck}), seems to give the most general SLLN to handle dependent
data. Since we will use such results to compare our findings to, we recall
one of the most sophisticated forms of the H\'{a}jek-R\'{e}nyi setting given
by T\'{o}m\'{a}cs and L\'{\i}bor (see \cite{thomas06}) denoted by $(GCHR)$.
These authors introduced a H\'{a}jek-R\'{e}nyi's inequality for
probabilities and, subsequently, got from it SLLN's for random sequences.
They obtained first :

\begin{theorem}
\label{theo2} Let $r$ be a positive real number, \ $a_{n}$ \ be a sequence
of nonnegative real numbers. Then the following two statements are
equivalent. \newline
(i) There exists \ $C>0$ \ such that for any \ ${n}\in \mathbb{\mathbf{%
\mathbb{N}
}}$ \ and any \ $\varepsilon >0$ 
\begin{equation*}
\mathbb{P}\left( {\max_{{\ell }\leq {n}}}|S_{\ell }|\geq \varepsilon \right)
\leq C\varepsilon ^{-r}{\sum_{{\ell }\leq {n}}}a_{\ell }.
\end{equation*}%
(ii) There exists $C>0$ such that for any nondecreasing sequence $%
(b_{n})_{n\in \mathbf{N}}$ of positive real numbers, for any ${n}\in \mathbb{%
\mathbf{%
\mathbb{N}
}}$ and any \ $\varepsilon >0$ 
\begin{equation*}
\mathbb{P}\left( {\max_{{\ell }\leq {n}}}|S_{\ell }|b_{\ell }^{-1}\geq
\varepsilon \right) \leq C\varepsilon ^{-r}{\sum_{{\ell }\leq {n}}}a_{\ell
}b_{\ell }^{-r}
\end{equation*}%
where $S_{n}=\sum_{i=1}^{n}X_{i}$ for all $n\in 
\mathbb{N}
.$
\end{theorem}

\noindent And next, they derived this SLLN from it.

\begin{theorem}
\label{theo3} Let $a_{n}\ $and$\ b_{n}$ be non-negative sequences of real
numbers and let $r>0$. Suppose that $b_{n}$ is a positive non-decreasing,
unbounded sequence of positive real numbers. Let us assume that%
\begin{equation*}
\sum_{n}\frac{a_{n}}{b_{n}^{r}}<+\infty
\end{equation*}%
and there exists $C>0$ such that for any ${n}\in \mathbb{\mathbf{%
\mathbb{N}
}}$ and any $\varepsilon >0$ 
\begin{equation*}
\mathbb{P}\left( \max_{{m}\leq {n}}|S_{m}|\geq \varepsilon \right) \leq C\
\varepsilon ^{-r}\sum_{{m}\leq {n}}a_{m}.
\end{equation*}%
Then 
\begin{equation*}
\lim_{{n}\rightarrow +\infty }\frac{S_{n}}{b_{n}}=0\ \ a.s.
\end{equation*}
\end{theorem}

\noindent For convenience, introduce these three notations. We say that a
sequence of random variables $X_{1},$ $X_{2},...$ has the $\mathbb{P}%
-max-variance(r)$ property, with $r>0,$ if and only if there exists a
constant $C>0$ such that for any $\ $fixed $n\geq 1,$ for any $\lambda >0,$ 
\begin{equation*}
\mathbb{P}\left( \max \left( |S_{1}|,...,|S_{n}|\right) \geq \lambda \right)
\leq C\text{ }\lambda ^{-r}\text{ }\mathbb{V}ar(S_{n}).
\end{equation*}%
It has the $\mathbb{V}ar-max-variance(r)$ property, $\ $with $r>0,$ if and
only if there exists a constant $C>0$ such that for any $\ $fixed $n\geq 1,$ 
\begin{equation*}
\mathbb{V}ar(\max (|S_{1}|,...,|S_{n}|))^{2/r}\leq C\text{ }\mathbb{V}%
ar(S_{n})
\end{equation*}%
and it has the $\mathbb{E}-max-variance(r)$ property,$\ $with $r>0,$ if and
only if there exists a constant $C>0$ such that for any $\ $fixed $n\geq 1,$ 
\begin{equation*}
\left( \mathbb{E}\left( \max (|S_{1}|,...,|S_{n}|)\right) ^{2}\right)
^{2/r}\leq C\text{ }\mathbb{V}ar(S_{n}).
\end{equation*}%
In the sequel we will say that $max-variance$ property is satisfied if one
of the three above $max-variance$ properties holds.

\bigskip

\noindent Theorem \ref{theo2} leads to these general laws.

\begin{proposition}
\label{prop1}Let $X_{1},$ $X_{2},...$ \ be a sequence of centered random
variables. Let $(b_{k})_{k\geq 1}$ be an increasing and nonbounded sequence
of positive real numbers. Assume that%
\begin{equation}
\lim \sup_{n\rightarrow +\infty }\sum_{1\leq i\leq n}b_{i}^{-r}cov(X_{i},%
\text{ }S_{n})<+\infty  \label{GCHR1}
\end{equation}%
and the sequence has the $\mathbb{P}-max-variance(r)$ property, $r>0$. Then $%
S_{n}/b_{n}\rightarrow 0$ $a.s.$ as $n\rightarrow +\infty .$ \bigskip

\noindent If the sequence has the $\mathbb{V}ar-max-variance(2)$ property or
the $\mathbb{E}-max-variance(2)$ property and if $\sum_{i\geq
1}b_{i}^{-2}\sum_{j\geq 1}Cov(X_{i},X_{j})<+\infty $, then $%
S_{n}/b_{n}\rightarrow 0$ $a.s.$ as $n\rightarrow +\infty .$
\end{proposition}

\bigskip

\begin{remark}
Here, (\ref{GCHR1}) is called the general condition of H\'{a}jek-R\'{e}nyi $%
(GCHR).$
\end{remark}

\bigskip \noindent \textbf{Proof.} If the sequence has the $\mathbb{E}%
-max-variance(r)$ property, then there exists a constant $C>0$ such that for
any $\ $fixed $n\geq 1,$ for any $\lambda >0$, and for $r=2,$%
\begin{equation*}
\mathbb{P}(\max (|S_{1}|,...,|S_{n}|)\geq \lambda )\leq \lambda ^{-r}\mathbb{%
V}ar(\max (|S_{1}|,...,|S_{n}|))
\end{equation*}%
\begin{equation*}
\leq \lambda ^{-r}\mathbb{E}\left( \max (|S_{1}|,...,|S_{n}|)\right) ^{2}
\end{equation*}%
\begin{equation*}
\leq C\lambda ^{-r}\mathbb{V}ar(S_{n})=C\lambda ^{-r}\sum_{i=1}^{n}\left[
\sum_{j=1}^{n}Cov(X_{i},\text{ }X_{j})\right] .
\end{equation*}%
The conclusion comes out by taking $a_{i}=\left[ \sum_{j=1}^{n}Cov(X_{i},%
\text{ }X_{j})\right] =Cov(X_{i},$ $S_{n})$ in the H\'{a}jek-R\'{e}nyi's
Theorem \ref{theo2} and applying Theorem \ref{theo3}.

\bigskip \noindent It is worth mentioning that the H\'{a}jek-R\'{e}nyi's
inequality is indeed very powerfull but, unfortunately, it works only if we
have the $max-variance$ property. For example, \ the $\mathbb{E}-max$
property holds for strictly stationary and associated sequences (see \cite%
{newman81}).

\bigskip

\noindent As to the squared indice method, it seems that it has not been
sufficiently standed to provide general strong laws for dependent data. We
aim at filling such a gap.

\bigskip

\noindent Indeed, in the next section, we provide a new general SLLN that
inspired by the squared indice method. This SLLN will be showed to have
interesting applications when comparing to the results of the present
section.

\section{Our results}

\noindent \label{sec1} In this section, we present a general SLLN based on
the squared indice method and give different forms in specific types of
dependent data including association with comparison with available results.
The result will be used in section \ref{sec2} to establish the strong
convergence for the continuous Hill's estimator with in the frame of EVT.

\begin{theorem}
\label{theo1}Let $X_{1},X_{2},\cdots $be an arbitrary sequence of rv's, and
let $\left( f_{i,n}\right) _{i\geq 1}$ be a sequence of measurable functions
such that $\mathbb{V}ar[f(X_{i,n})]<+\infty $, for $i\geq 1$ and $n\geq 1.$
If for some $\delta ,0<\delta <3,$ 
\begin{equation}
C_{1}=\sup_{n\geq 0}\sup_{q\geq 1}\mathbb{V}ar\left( \dfrac{1}{q^{(3-\delta
)/4}}\sum_{i=1}^{q}f_{i,\text{ }n}(X_{i})\right) <+\infty  \label{GCIP1}
\end{equation}%
and for some $\delta ,$ $0<\delta <3,$ 
\begin{equation}
C_{2}=\sup_{n>0}\sup_{k\geq 1}\sup_{q\text{ : }q^{2}+1\leq k\leq
(q+1)^{2}}\sup_{1\leq j\leq k}\mathbb{V}ar\left( \dfrac{1}{q^{(3-\delta )/2}}%
\sum_{i=1}^{j-q^{2}+1}f_{q^{2}+i,\text{ }n}\left( X_{q^{2}+i}\right) \right)
<+\infty  \label{GCIP2}
\end{equation}%
hold, then 
\begin{equation*}
\dfrac{1}{n}\sum_{i=1}^{n}\left( f_{i,\text{ }n}(X_{i})-\mathbb{E}\left(
f_{i,\text{ }n}(X_{i})\right) \right) \longrightarrow
0\;a.s\;as\;n\longrightarrow +\infty .
\end{equation*}
\end{theorem}

\begin{remark}
\bigskip We say that the sequence $X_{1},X_{2},\cdots,X_{n}$ satisfies the $%
(GCIP)$ whenever $(\ref{GCIP1})$ and $(\ref{GCIP2})$ hold.
\end{remark}

\noindent \textbf{Proof.} It suffices to prove the announced results for $%
Y_{i}=f_{i,\text{ }n}(X_{i})$ and $\mathbb{E}(Y_{i})=0,$ $i\geq 1$. Observe
that omitting the subscript $n$ does not cause any ambiguity in the proof
below. We have 
\begin{equation*}
\mathbb{P}\left( \left\vert \dfrac{1}{k}\sum_{i=1}^{k}Y_{i}\right\vert
>k^{-\beta }\right) \leq \mathbb{P}\left( \left\vert
\sum_{i=1}^{k}Y_{i}\right\vert >k^{1-\beta }\right) \leq \frac{1}{%
k^{2(1-\beta )}}\mathbb{V}ar\left( \sum_{i=1}^{k}Y_{i}\right) .
\end{equation*}%
\noindent We apply this formula for $k=q^{2}$ and get for $0<\delta <3,$

\begin{equation*}
\mathbb{P}\left( \left\vert \frac{1}{q^{2}}\sum_{i=1}^{q^{2}}Y_{i}\right%
\vert >q^{-2\beta }\right) \leq \frac{1}{q^{4(1-\beta )}}\mathbb{V}ar\left(
\sum_{i=1}^{q^{2}}Y_{i}\right)
\end{equation*}%
\begin{equation*}
\leq \dfrac{1}{q^{1+\delta -4\beta }}\mathbb{V}ar\left( \dfrac{1}{%
q^{(3-\delta )/2}}\sum_{i=1}^{q^{2}}Y_{i}\right) ,
\end{equation*}%
and there exists $q_{_{0}}\in 
\mathbb{N}
^{\ast }$ such that for $q>q_{_{0}}$ 
\begin{equation*}
\mathbb{P}\left( \left\vert \frac{1}{q^{2}}\sum_{i=1}^{q^{2}}Y_{i}\right%
\vert >q^{-2\beta }\right) <\dfrac{C_{1}}{q^{1+\delta -4\beta }}.
\end{equation*}%
Then we have for $\beta <(\delta /4)$, $\sum_{q=1}^{+\infty }\mathbb{P}%
\left( \left\vert \frac{1}{q^{2}}\sum_{i=1}^{q^{2}}Y_{i}\right\vert
>q^{-2\beta }\right) <+\infty $. We conclude that 
\begin{equation}
\frac{1}{q^{2}}\sum_{i=1}^{q^{2}}Y_{i}\longrightarrow 0\;a.s.  \label{f1}
\end{equation}%
Now set $q^{2}\leq k\leq (q+1)^{2}$ and $\epsilon _{k,\text{ }q}=0$ if $%
k=q^{2}$ and $1$ otherwise. We have

\begin{equation*}
\frac{1}{k}\sum_{i=1}^{k}Y_{i}-\frac{1}{q^{2}}\sum_{i=1}^{q^{2}}Y_{i}=\frac{1%
}{k}\sum_{i=1}^{k}Y_{i}-\frac{1}{k}\sum_{i=1}^{q^{2}}Y_{i}+\frac{1}{k}%
\sum_{i=1}^{q^{2}}Y_{i}-\frac{1}{q^{2}}\sum_{i=1}^{q^{2}}Y_{i}
\end{equation*}

\begin{equation*}
=\dfrac{\epsilon _{k,q}}{k}\left(
\sum_{i=1}^{k}Y_{i}-\sum_{i=1}^{q^{2}}Y_{i}\right) +\dfrac{1}{q^{2}}%
\sum_{i=1}^{q^{2}}Y_{i}\left( \dfrac{q^{2}-k}{k}\right)
\end{equation*}%
\begin{equation}
=\dfrac{\epsilon _{k,\text{ }q}}{k}\left( \sum_{i=q^{2}+1}^{k}Y_{i}\right) +%
\dfrac{1}{q^{2}}\sum_{i=1}^{q^{2}}Y_{i}\left( \dfrac{q^{2}-k}{k}\right).
\label{f2}
\end{equation}%
But $(q^{2}-k)/k\rightarrow 0$ as $q\rightarrow +\infty $. This combined
with (\ref{f1}) proves that the second term of (\ref{f2}) converges to zero $%
a.s.$ It remains to handle the first term. Notice that we let $q\rightarrow
+\infty $ and then $k(q)$ defined by $q(k)^{2}\leq k\leq (q(k)+1)^{2}$ also
goes to infinity. Then, for $0<\delta <3,$

\begin{equation*}
\mathbb{P}\left( \frac{1}{k}\left\vert \epsilon _{k,\text{ }%
q}\sum_{i=q^{2}+1}^{k}Y_{i}\right\vert >k^{1-\beta }\right) \leq \mathbb{P}%
\left( \left\vert \epsilon _{k,\text{ }q}\sum_{i=q^{2}+1}^{k}Y_{i}\right%
\vert >k^{1-\beta }\right)
\end{equation*}%
\begin{equation*}
\leq \mathbb{P}\left( \left\vert \epsilon _{k,\text{ }q}%
\sum_{i=q^{2}+1}^{k}Y_{i}\right\vert >q^{2(1-\beta )}\right) \leq \frac{%
\epsilon _{k,\text{ }q}}{q^{4-4\beta }}\mathbb{V}ar\left(
\sum_{i=q^{2}+1}^{k}Y_{i}\right)
\end{equation*}%
\begin{equation*}
\leq \frac{\epsilon _{k,\text{ }q}}{q^{1+\delta -4\beta }}\mathbb{V}ar\left( 
\frac{1}{q^{(3-\delta )/2}}\sum_{i=q^{2}+1}^{k}Y_{i}\right) \leq \frac{%
\epsilon _{k,\text{ }q}C_{2}}{q^{1+\delta -4\beta }}.
\end{equation*}%
Now for $0<\beta <\delta /4,\sum_{k=1}^{+\infty }\mathbb{P}\left( \epsilon
_{k,\text{ }q}\left\vert \sum_{i=q^{2}+1}^{k}Y_{i}\right\vert >k^{1-\beta
}\right) <+\infty .$ Then 
\begin{equation}
\frac{\epsilon _{k,\text{ }q}}{k}\left[ \sum_{i=1}^{k}Y_{i}-%
\sum_{i=1}^{q^{2}}Y_{i}\right] \longrightarrow 0\;a.s.  \label{f3}
\end{equation}%
Now in view of $(\ref{f1})$, $(\ref{f2})$ and $(\ref{f3})$ and since ($%
q^{2}-k)/k\rightarrow 0$, we may conclude the proof.

\begin{remark}
In most cases, conditions $(\ref{GCIP1})$ and $(\ref{GCIP2})$ are used for $%
\delta =1$, as it is the case for the independent and indentically
distributed random variables. We will exhibit a situation in Proposition \ref%
{prop2} that cannot be handled without using $(\ref{GCIP1})$ and $(\ref%
{GCIP2})$ for $\delta <1.$
\end{remark}

\subsection{Comparison and particular cases}

\label{subsec4} Let us see how $(GCIP)$, that is fulfilment of conditions $(%
\ref{GCIP1})$ and $(\ref{GCIP2})$, works in special cases. We have to
compare our $(GCIP)$ to $(GCHR).$ But $(GCHR)$ is used only when $%
max-variance$ property is satisfied. We only consider the case where $%
X_{1},X_{2},\cdots $ are real and the $f_{i}$'s are identity functions.

\subsubsection{Independence case.}

\label{subsec5} By using Theorem \ref{theo2}, we observe that we have the $%
\mathbb{P}-max-variance(2)$ property, that is the Kolmogorov's maximal
inequality. By using the H\'{a}jek-R\'{e}nyi's general condition, we have
the strong law of large numbers of Kolmogorov : $S_{n}/n\rightarrow 0$ $a.s.$
whenever%
\begin{equation}
\sum_{n\geq 1}\mathbb{V}ar(X_{n})/n^{2}<+\infty .  \label{CI1}
\end{equation}%
To apply Theorem \ref{theo1} here, we notice that the sequence of variances $%
\mathbb{V}ar(S_{n})$ is non-decreasing in $n.$ Then $(\ref{GCIP1})$ and $(%
\ref{GCIP2})$ are implied by, for some $0<\nu _{1}$ and $0<\nu _{2},$%
\begin{equation*}
\sup_{k\geq 1}\frac{1}{k^{1+\nu _{1}}}\sum_{i=1}^{k}\mathbb{V}%
ar(X_{i})<+\infty \text{ and }\sup_{k\geq 1}\frac{1}{k^{2+\nu _{2}}}%
\sum_{i=k^{2}+1}^{(k+1)^{2}}\mathbb{V}ar(X_{i})<+\infty .
\end{equation*}%
But, by observing that the latter is 
\begin{equation*}
k^{-(2+\nu _{2})}\sum_{i=k^{2}}^{(k+1)^{2}}\mathbb{V}ar(X_{i})=k^{-(2+\nu
_{2})}\left[ \sum_{i=1}^{(k+1)^{2}}\mathbb{V}ar(X_{i})-\sum_{i=1}^{k^{2}}%
\mathbb{V}ar(X_{i})\right] ,
\end{equation*}%
we conclude that the SLLN is implied by%
\begin{equation}
\sup_{k\geq 1}\frac{1}{k^{1+\nu }}\sum_{i=1}^{k}\mathbb{V}ar(X_{i})<+\infty ,
\label{CI2}
\end{equation}%
for some $\nu >0.$ In the independant case, one has the SLLN for $%
k^{-1}\sum_{i=1}^{k}\mathbb{V}ar(X_{i})\rightarrow \sigma ^{2}.$ And the
parameter $\nu $ in $(\ref{CI2})$ is useless in that case. But\ the
availability of the parameter $\nu $ is important for situations beyond the
classical cases. As a first example, let us use the Kolmogorov's Theorem and
construct a probability space holding a sequence of independent centered
rv's $X_{1},X_{2},\cdots $ with $\mathbb{E}X_{n}^{2}=n^{1/3}.$ But $(\ref%
{CI2})$ does not hold for $\nu =0$ since%
\begin{equation*}
\frac{1}{n}\sum_{i=1}^{n}i^{1/3}\geq \frac{1}{n}\int_{1}^{n}x^{1/3}dx\geq 
\frac{3}{4}\left( n^{1/3}-1\right) \rightarrow +\infty ,\text{ as }%
n\rightarrow +\infty
\end{equation*}%
while $(GCHR)$ entails the SLLN.

\bigskip \noindent We will consider in proposition \ref{prop2} below an
important other example which cannot be concluded unless we use a positive
value of $\nu .$ Now, if we may take $\nu =1/3,$ we have that $n^{-\left(
1+\nu \right) }\sum_{i=1}^{n}i^{1/3}$ is bounded and our Theorem also
ensures the SLLN. \bigskip

\noindent Now if the sequence is second order stationary, then $(\ref{GCIP1}%
) $ and $(\ref{GCIP2})$ are both valid. Also, if the variances are bounded
by a common constant $C_{0}$, both $(\ref{GCIP1})$ and $(\ref{GCIP2})$ are
valid.

\subsubsection{Pairwise negatively dependent variables.}

\label{subsec6} In that case, we may drop the covariances in $(GCIP)$ and
then (\ref{GCIP1}) and (\ref{GCIP2}) lead to $(\ref{CI2})$ as a general
condition for the validity of the SLLN in the independent case. As to $%
(GCHR) $, we don't have any information whether or not the $max-variance$
property holds.

\subsubsection{Associated sequences}

\label{subsec7} Here $\mathbb{V}ar(S_{n})$ is non-decreasing in $n$ and $%
(GCIP)$ becomes for $\nu =(1-\delta )/2\geq 0$ with $0<\delta <1$%
\begin{equation}
\sup_{q\geq 1}\frac{1}{q^{1+\nu }}\mathbb{V}ar\left(
\sum_{i=1}^{q}X_{i}\right) <+\infty  \label{gca1}
\end{equation}%
and%
\begin{equation}
\sup_{q\geq 1}\frac{1}{q^{2(1+\nu )}}\mathbb{V}ar\left(
\sum_{i=q^{2}+1}^{(q+1)^{2}}X_{i}\right) <+\infty .  \label{gcal2}
\end{equation}%
If the sequence is second order stationary, then (\ref{gca1}) implies (\ref%
{gcal2}), since%
\begin{equation*}
\frac{1}{q^{2(1+\nu )}}\mathbb{V}ar\left(
\sum_{i=q^{2}+1}^{(q+1)^{2}}X_{i}\right) =\frac{\left( 2q+1\right) ^{1+\nu }%
}{q^{2(1+\nu )}}\left[ \frac{1}{\left( 2q+1\right) ^{1+\nu }}\mathbb{V}%
ar\left( \sum_{i=1}^{2q+1}X_{i}\right) \right]
\end{equation*}%
\begin{equation*}
\sim \frac{2}{q^{(1+\nu )}}\mathbb{V}ar\left( \frac{1}{k^{\left( 1+\nu
\right) /2}}\sum_{i=1}^{k}X_{i}\right) ,
\end{equation*}%
for $k=2q+1.$ And (\ref{gca1}) may be witten as 
\begin{equation}
\sup_{q\geq 1}\frac{1}{q^{\nu }}\left[ \mathbb{V}ar(X_{1})+\frac{2}{q}%
\sum_{i=2}^{q}\left( q-i+1\right) Cov(X_{1},\text{ }X_{i})\right] <+\infty .
\label{q2}
\end{equation}%
This is our general condition under which SLLN holds for second order
stationary associated sequence. Then, by the Kronecker lemma, we have the
SLLN if 
\begin{equation}
\sigma ^{2}=\mathbb{V}ar(X_{1})+2\sum_{i=2}^{+\infty }Cov(X_{1},\text{ }%
X_{i})<+\infty .  \label{new01}
\end{equation}%
Condition (\ref{new01}) is obtained by Newman \cite{newman81}. Clearly, by
the Ces\`{a}ro lemma, (\ref{new01}) implies%
\begin{equation}
\lim_{q\rightarrow +\infty }\frac{1}{q}\sum_{i=1}^{q}Cov(X_{1},\text{ }%
X_{i})\rightarrow 0.  \label{newmanlebowitz}
\end{equation}%
And, in fact, the latter is a necessary condition of strong law of large
numbers as proved in Theorem 7 in \cite{newman84}, from the original result
of Lebowitz (see \cite{lebow}).

\bigskip

\noindent The reader may find a larger review on this subject in \cite{rao}.
Our result seems more powerful since we may still have the strong law of
large numbers even if $\sigma ^{2}=+\infty $.

\bigskip

\noindent We only need to check condition (\ref{q2}). We will comment this
again after Proposition \ref{prop2}.

\bigskip

\noindent For strictly stationary associated sequences with finite variance,
we have the $\mathbb{E}-max-variance(2)$ property (see \cite{newman81}).
Then $(GCHR)$ may be used. It becomes%
\begin{equation}
\lim \sup_{n}\sum_{i=1}^{n}\frac{1}{i^{2}}Cov(X_{i},\text{ }S_{n})<+\infty ,
\label{cgchr}
\end{equation}%
which is equivalent to 
\begin{equation*}
\lim \sup_{n}\left[ \sum_{i=1}^{n}\frac{\mathbb{V}ar(X_{i})}{i^{2}}%
+\sum_{j=2}^{n}\left( \sum_{i=1}^{n-j+1}\frac{1}{i^{2}}+\sum_{i=j}^{n}\frac{1%
}{i^{2}}\right) Cov(X_{1},\text{ }X_{j})\right] <+\infty
\end{equation*}%
and reduces to%
\begin{equation*}
\sum_{j=2}^{+\infty }Cov(X_{1},\text{ }X_{j})<+\infty .
\end{equation*}%
We then see that $(GCHR)$ gives weaker results than ours. Indeed, in our
formula (\ref{q2}), we did not require that $\frac{2}{q}\sum_{i=2}^{q}\left(
q-i+1\right) Cov(X_{1},$ $X_{i})$ is bounded. It may be allowed to go to
infinity at a slower convergence rate than $q^{-\nu }$. Then our condition (%
\ref{q2}) besides being more general, applies to any associated sequences
and is significantly better than the $(GCHR)$ for strictly stationary
sequences.

\bigskip

\noindent Nevertheless, for (\ref{cgchr}), it is itself more powerfull than
Theorem 6.3.6 and Corollary 6.3.7 in \cite{rao}, due to the use of Theorem %
\ref{theo2} and Proposition \ref{prop1}, \ of T\'{o}m\'{a}cs and L\'{\i}bor $%
\left( \text{see \cite{thomas06}}\right) $. Such a result is also obtained
by Yu (1993) $\left( \text{see \cite{yu93}}\right) $ for the strong
convergence of empirical distribution function for associated sequence with
identical and continuous distribution.

\bigskip

\noindent Birkel (see \cite{birkel}) used direct computations on the
convariance structure for associated variables and got the following
condition%
\begin{equation*}
\lim \sup_{n}\sum_{i=1}^{n}\frac{1}{i^{2}}Cov(X_{i},\text{ }S_{i})<+\infty
\end{equation*}%
for SLLN for associated variables.

\bigskip

\noindent Now, to sum up, the comparison between $(GCIP)$ and $(GCHR)$ is as
follows

\bigskip

\begin{itemize}
\item[1] For independent case the two conditions are equivalent.\newline

\item[2] In negatively associated case, the form of $(GCIP)$ for independent
case remains valid. And we have no information whether the $max-variance$
property holds to be able to apply $(GCHR).$\newline

\item[3] For association with strictly stationary of sequences, $(GCIP)$
gives a better condition than $(GCHR)$.\newline

\item[4] For association with no information on stationarity, so $(GCHR)$
cannot be applied unless a $max-variance$ property is proved. Our condition
still works and is the same as for the stationary associated sequences in
point 3.

\item[5] For arbitrary sequences with finite variances, point 4 may be
recontacted.
\end{itemize}

\bigskip

\noindent In conclusion, our method effectively brings a significant
contribution to SLLN for associated random variables. And we are going to
apply it to an associated sequence in the extreme value theory fields.

\section{\protect\bigskip Applications}

\label{sec2}

\subsection{Application to extreme value theory}

\label{subsec8} The EVT offers us the opportunuity to directly apply our
general conditions $(\ref{GCIP1})$ and $(\ref{GCIP2})$ to a sum of dependent
and non-stationary random variables and to show how to proceed in such a
case.

\bigskip

\noindent We already emphasized the importance of the parameter $\nu
=(1-\delta )/2$ in $(GCIP)$. In the example we are going to treat, we will
see that a conclusion cannot be achieved with $\nu =0.$

\bigskip

\noindent Let $E_{1},E_{2},...$ be an infinite sequence of independent
standard exponential random variables, $f(j)$ is an increasing function of
the integer $j\geq 0$ with $f(0)=0$ and $\gamma >0$ a real parameter. Define
the following sequences of random variables 
\begin{equation}
W_{k}=\sum_{j=1}^{k-1}f(j)\left[ \exp \left( -\gamma
\sum_{h=j+1}^{k-1}E_{h}/h\right) -\exp \left( -\gamma
\sum_{h=j}^{k-1}E_{h}/h\right) \right] ,\text{ }k\geq 1.  \label{sum01}
\end{equation}%
The characterization of the asymptotic behavior of (\ref{sum01}) has
important applications and consequences in two important fields : the
extreme value theory in statistics and the central limit theorem issue for
sum of non stationary associated random variables. Let us highlight each of
these points.

\bigskip

\noindent On one side, let $X,X_{1},X_{2},...$ \ be independent and
identically random variables in Weibull extremal domain of parameter $\gamma
>0$ such that $X>0$ and let $X_{1,n}\leq X_{2,n}\leq ...\leq X_{n,n}$ denote
the order statistics based on the $n\geq 1$ observations. The distribution
function $G$ of $Y=\log X$ has a finite upper endpoint $y_{0}$ and admits
the following representation :%
\begin{equation}
y_{0}-G^{-1}(1-u)=cu^{1/\gamma }(1+p(u))\exp \left(
\int_{u}^{1}t^{-1}b(t)dt\right) ,\text{ }u\in (0,1)  \label{rep01}
\end{equation}%
where $c$ is a constant and, $p(u)$ and $b(u)$ are functions of $u\in
(\bigskip 0.1)$ such that $(p(u),b(u))\rightarrow 0$ as $u\rightarrow 0.$
This is called a representation of a sequence of random variables in the
Weibull domain of attraction.

\bigskip

\noindent To stay simple, suppose that $p(u)=b(u)=0$ for all $u\in (0,1)$
consider the simplest case%
\begin{equation}
y_{0}-G^{-1}(1-u)=u^{\gamma },\text{ }u\in (0,1).  \label{rep02}
\end{equation}%
The so-called Hill's statistic, based on the identity function $id(x)=x$ and
the $k$ largest values with $1\leq k\leq n,$%
\begin{equation}
T_{n}(id)=\frac{1}{id(k)}\sum_{j=1}^{k}id(j)\left( \log X_{n-j+1,\text{ }%
n}-\log X_{n-j,\text{ }n}\right)  \label{hillsimp}
\end{equation}%
is an estimator of $\gamma $ in the sense that 
\begin{equation*}
\frac{T_{n}(id)}{(y_{0}-\log X_{n-k,\text{ }n})}\rightarrow ^{\mathbb{P}%
}(\gamma +1)^{-1},
\end{equation*}%
as $n\rightarrow +\infty .$ When we replace the identity function with an
increasing function $f(j)$ of the integer $j\geq 0$ with $f(0)=0,$ we get
the functional Hill's estimator defined as%
\begin{equation}
T_{n}(f)=\frac{1}{f(k)}\sum_{j=1}^{k}f(j)\left( \log X_{n-j+1,\text{ }%
n}-\log X_{n-j,\text{ }n}\right)  \label{hillfonct}
\end{equation}%
introduced by D\`{e}me E., Lo G.S. and Diop, A. (2012) (see \cite{demelo12}%
). From this processus is derived the Diop and Lo (2006) (see \cite{dioplo})
generalization of Hill's statistic. We are going to highlight that $%
f(k)T_{n}(f)/(y_{0}-\log X_{n-k,\text{ }n})$ is of the form of (\ref{sum01})
when (\ref{rep02}) holds. \ We have to use two representations. The R\'{e}%
nyi's representation allows to find independent standard uniform random
variables $U_{1},U_{2},...$ such that the following equalities in
distribution hold%
\begin{equation*}
\{\log Y_{j},\text{ }j\geq 1\}=_{d}\{G^{-1}(1-U_{j}),\text{ }j\geq 1\}
\end{equation*}%
and 
\begin{equation*}
\left\{ \{\log X_{n-j+1,\text{ }n},\text{ }1\leq j\leq n\},\text{ }n\geq
1\right\} =_{d}\left\{ \{G^{-1}(1-U_{j,\text{ }n}),\text{ }1\leq j\leq n\},%
\text{ }n\geq 1\right\} .
\end{equation*}%
Next, by the Malmquist representation (see (\cite{shwell}), p. 336), we have
for each $n\geq 1,$ the following equality in distribution holds%
\begin{equation*}
\{j^{-1}\log (U_{j+1,\text{ }n}/U_{j,\text{ }n}),\text{ }1\leq j\leq
n\}=_{d}\{E_{j}^{(n)},\text{ }1\leq j\leq n\},
\end{equation*}%
where $E_{j}^{(n)},$ $1\leq j\leq n,$ are independent exponential random
variables. We apply these two tools to get that for each fixed $n$ and $%
k=k(n)$%
\begin{equation}
\frac{T_{n}(f)}{(y_{0}-\log X_{n-k,\text{ }n})}=_{d}W_{k(n)}.  \label{ext}
\end{equation}%
For an arbitrary element of the Weibull extremal domain of attraction, it
may be easily showed that $f(k)T_{n}(f)/(y_{0}-\log X_{n-k,\text{ }n})$ also
behaves as (\ref{sum01}) if some extra conditions are imposed of the
auxiliary functions $p$ and $b.$ Hence a complete characterization of the
asymptotic bevahior of (\ref{sum01}) provides asymptotic laws in extreme
value theory\bigskip . \noindent On another side, easy algebra leads to 
\begin{equation*}
W_{k}=f(k-1)+\sum_{j=1}^{k-1}\Delta f(j)\exp \left( -\gamma
\sum_{h=j}^{k-1}E_{h}/h\right) ,
\end{equation*}%
where $\Delta f(j)=f(j)-f(j-1),$ $j\geq 1.$ We consider%
\begin{equation}
W_{k}^{\ast }=W_{k}-\mathbb{E}(W_{k})=\sum_{j=1}^{k-1}\Delta f(j)\left[ \exp
\left( -\gamma \sum_{h=j}^{k-1}E_{h}/h\right) -\mathbb{E}\exp \left( -\gamma
\sum_{h=j}^{k-1}E_{h}/h\right) \right] .  \label{sum03}
\end{equation}%
This is a sum of non stationary dependent random variables. In fact the $rv$%
's 
\begin{equation*}
\Delta f(j)\left[ \exp \left( -\gamma \sum_{h=j}^{k-1}E_{h}/h\right) -%
\mathbb{E}\exp \left( -\gamma \sum_{h=j}^{k-1}E_{h}/h\right) \right]
\end{equation*}%
are associated.

\bigskip

\noindent Now, we are going to apply our general conditions to ( \ref{sum04}%
), defined below 
\begin{equation}
S_{k}^{\ast }=\sum_{j=1}^{k-1}\Delta f(j)\left[ \exp \left( -\gamma
\sum_{h=j}^{k-1}E_{h}/h\right) -\mathbb{E}\exp \left( -\gamma
\sum_{h=j}^{k-1}E_{h}/h\right) \right] \alpha (k),  \label{sum04}
\end{equation}%
where $\alpha (k)$ is a sequence of positive real numbers. Next, we will
particularize the result for $f(j)=j^{\tau },$ $\tau >0.$ Our results depend
on computation techniques developed in \cite{adja}. Here are our results :

\begin{proposition}
\label{prop2} Suppose that, for $L$ and $q$ large enough such that $L\leq
q^{2},$ the following conditions hold for some $\delta ,$ $0<\delta <3.$%
\begin{equation}
\sup_{k\geq L}\frac{\alpha ^{2}(k)}{k^{2\gamma +1+\nu }}\sum_{j=L}^{k-1}%
\Delta ^{2}f(j)j^{2\gamma }<+\infty ,  \label{1}
\end{equation}%
\begin{equation}
\sup_{k\geq L}\frac{\alpha ^{2}(k)}{k^{1+\nu }}\sum_{j=L+1}^{k-1}\left[
\sum_{i=L}^{j-1}\Delta f(i)\right] \Delta f(j)\frac{1}{j}<+\infty ,
\label{2}
\end{equation}%
\begin{equation}
\sup_{k\geq L}\frac{\alpha ^{2}(k)}{k^{1+\nu }}\sum_{L\leq j\leq k-1}\Delta
f(j)/j<+\infty ,  \label{2m}
\end{equation}%
\begin{equation}
\sup_{k\geq 1}\frac{1}{q^{(3-\delta )}}\sum_{i=1}^{2q+1}\alpha ^{2}(k)\Delta
^{2}f(q^{2}+i)\left( \frac{q^{2}+i}{k}\right) ^{2\gamma }<+\infty  \label{3}
\end{equation}%
and%
\begin{equation}
\sup_{k\geq 1}\sup_{(q^{2}+1)\leq k\leq (q+1)^{2}}\frac{\alpha ^{2}(k)}{%
q^{(3-\delta )}}\sum_{j=2}^{2q+1}\left[ \sum_{i=1}^{j-1}\Delta f(q^{2}+i)%
\right] \Delta f(q^{2}+j)\frac{1}{q^{2}+j}<+\infty .  \label{4}
\end{equation}%
Then 
\begin{equation*}
\frac{S_{k}^{\ast }}{k}\rightarrow 0\text{ }a.s.
\end{equation*}%
Further, if \ 
\begin{equation*}
\mu _{k}=\sum_{j=1}^{k-1}\alpha (k)\Delta f\left( j\right) \mathbb{E}\exp
\left( -\gamma \sum_{h=j}^{k-1}E_{h}/h\right) \rightarrow \mu ,
\end{equation*}%
where $\mu $ is a finite$,$ then%
\begin{equation*}
k^{-1}\sum_{j=1}^{k-1}\alpha (k)\Delta f\left( j\right) \exp \left( -\gamma
\sum_{h=j}^{k-1}E_{h}/h\right) \rightarrow \mu \text{ }a.s.
\end{equation*}
\end{proposition}

\begin{proposition}
\label{prop3} For $f(j)=j^{\tau },$ if (\ref{1}), (\ref{2}), (\ref{2m}), (%
\ref{3}) and (\ref{4}) hold, $\alpha (k)=1/k^{\tau -1}$ and if $\mu =\tau
/(\tau +\gamma ).$ Then%
\begin{equation*}
\frac{1}{k^{\tau }}\sum_{j=1}^{k-1}\left( j^{\tau }-(j-1)^{\tau }\right)
\exp \left( -\gamma \sum_{h=j}^{k-1}E_{h}/h\right) \rightarrow \frac{\tau }{%
\gamma +1}\text{ }a.s.\text{ as }k\rightarrow +\infty .
\end{equation*}
\end{proposition}

\begin{remark}
\bigskip Since these results are only based on moments, the a.s. convergence
remains true for $T_{n}(f)/(y_{0}-\log X_{n-k,n})$ in vertue of (\ref{ext}).
We get under the model \ that%
\begin{equation*}
\frac{T_{n}(f)}{(y_{0}-\log X_{n-k,n})}\rightarrow \frac{\tau }{\gamma +1}%
\text{ }a.s.\text{ as }n\rightarrow +\infty \text{ and }k=k(n)\rightarrow
+\infty \text{ and }k/n\rightarrow 0
\end{equation*}%
under the assumptions (\ref{1}), (\ref{2}), (\ref{2m}), (\ref{3}) and (\ref%
{4}), in the general case.
\end{remark}

\begin{remark}
This strong law may be easily checked by Monte Carlo simulations. For
example, consider $\gamma =2$ and $\tau =1.$ We observe the following errors
corresponding to the values of $50$, $75$ and $100$ of $k$ : $0.358$, $0.321$
and $0.3332$. This shows the good performance of this strong law for the
particular values $\gamma =2$ and $\tau =1.$
\end{remark}

\subsubsection{Proofs}

Both proofs of the two propositions are postponed in the Appendix Section .

\section{Conclusion and Perspectives}

\label{conclu} \noindent We have established a general SLLN and applied it
to associated variables. Comparison with SLLN's derived from the H\'{a}jek-R%
\'{e}nyi inequality proved that this SLLN is not trivial. We have also used
it to find the strong convergence of statistical estimators under
non-stationary associated samples in EVT.

\bigskip \noindent It seems that it has promising applications in
non-parametric statistic, when dealing with the strong convergence of the
empirical process and the non-parametric density estimator for a stationary
sequence with an arbitrary parent distribution function.

\section{Appendix}

\label{sec3}

\subsection{Proofs of Propositions \protect\ref{prop2} and \protect\ref%
{prop3}}

\subsubsection{ Assumptions}

We have to show that the assumptions of Proposition \ref{prop2} entail the
general condition $(GCIP)$. We first remind that 
\begin{equation*}
S_{k}^{\ast }=\sum_{j=1}^{k-1}\Delta f(j)\left[ \exp \left( -\gamma
\sum_{h=j}^{k-1}E_{h}/h\right) -\mathbb{E}\exp \left( -\gamma
\sum_{h=j}^{k-1}E_{h}/h\right) \right] \alpha (k)
\end{equation*}%
\bigskip that we write as 
\begin{equation}
S_{k}^{\ast }=\sum_{j=1}^{k-1}\alpha (k)\Delta f(j)\left(
S_{j,k}-s_{j,k}\right) ,  \label{eq}
\end{equation}%
where $S_{j,\text{ }k}=\exp \left( -\gamma \sum_{h=j}^{k-1}E_{h}/h\right) $
and $s_{j,\text{ }k}=\mathbb{E}\exp \left( -\gamma
\sum_{h=j}^{k-1}E_{h}/h\right) $. Next, we are going to check (\ref{GCIP1})
and (\ref{GCIP2}) for this sum of random variables. Fix $\delta ,$ $0<\delta
<3.$ Let us split (\ref{eq}) into%
\begin{equation*}
S_{k}^{\ast }=\sum_{j=1}^{L-1}\alpha (k)\Delta f(j)\left( S_{j,\text{ }%
k}-s_{j,\text{ }k}\right) +\sum_{j=L}^{k-1}\alpha (k)\Delta f(j)\left( S_{j,%
\text{ }k}-s_{j,\text{ }k}\right)
\end{equation*}%
\begin{equation*}
=:S_{L}^{1}+S_{L}^{2}.
\end{equation*}%
Then for $\nu =(1-\delta )/2$ with $0<\delta <1,$ 
\begin{equation*}
\frac{1}{k^{1+\nu }}\mathbb{V}ar(S_{k}^{\ast })=\frac{1}{k^{1+\nu }}\mathbb{V%
}ar\left( S_{L}^{1}\right) +\frac{1}{k^{1+\nu }}\mathbb{V}ar\left(
S_{L}^{2}\right) +\frac{2}{k^{1+\nu }}cov\left( S_{L}^{1},S_{L}^{2}\right)
\end{equation*}%
\begin{equation*}
=:A_{k}+B_{k}+2C_{k}.
\end{equation*}%
Let us treat each term in the above equality. Here, we use Formulas 18 and
21\ \ in \cite{adja} and take $L$ large enough to ensure%
\begin{equation}
\mathbb{V}ar\left( S_{j,\text{ }k}\right) =\left( \frac{j}{k-1}\right)
^{2\gamma }V(1,j)V(2,j),  \label{fadja01}
\end{equation}%
with%
\begin{equation*}
\left\vert V(1,j)\right\vert =1+O(j^{-1})\text{ and }0\leq V(2,j)\leq \frac{%
2\gamma ^{2}\left\vert a_{1}(\in )\right\vert }{j}.
\end{equation*}%
And%
\begin{equation}
Cov(S_{j,\text{ }k},S_{j+\ell ,\text{ }k})=\mathbb{V}ar\left( S_{j+\ell ,%
\text{ }k}\right) \left( \frac{j}{j+\ell -1}\right) ^{\gamma }(1+O(j^{-1})).
\label{fadja02}
\end{equation}%
We suppose that $L$ is large enough so that $\left\vert V(1,j)\right\vert
\leq 1/2,$ for $j\geq L.$ \bigskip First we see that 
\begin{equation}
A_{k}\rightarrow 0,\text{ }as\text{ }k\rightarrow +\infty ,  \label{gc01}
\end{equation}%
since $\mathbb{V}ar(S_{L}^{1})$ is let constant with $L.$ Next, split $B_{k}$
into 
\begin{equation*}
B_{k}=\frac{1}{k^{1+\nu }}\sum_{j=L}^{k-1}\alpha ^{2}(k)\Delta ^{2}f(j)%
\mathbb{V}ar\left( S_{j,\text{ }k}-s_{j,\text{ }k}\right)
\end{equation*}%
\begin{equation*}
+\frac{1}{k^{1+\nu }}\sum_{L\leq i\neq j\leq k-1}\alpha ^{2}(k)\Delta
f(j)\Delta f(i)Cov\left( S_{i,\text{ }k},S_{j,\text{ }k}\right)
\end{equation*}%
\begin{equation*}
=:B_{k,\text{ }1}+B_{k,\text{ }2}.
\end{equation*}%
By (\ref{fadja01}) we get

\begin{equation}
B_{k,\text{ }1}=\frac{1}{k^{1+\nu }}\sum_{j=L}^{k-1}\alpha ^{2}(k)\Delta
^{2}f(j)\mathbb{V}ar\left( S_{j,\text{ }k}\right) \leq (1/2)\frac{\alpha
^{2}(k)}{k^{2\gamma +1+\nu }}\sum_{j=L}^{k-1}\Delta ^{2}f(j)j^{2\gamma }.
\label{gc02}
\end{equation}%
\noindent Now let us turn to the term $B_{k,\text{ }2}.$ Let us remark that
the rv$^{\prime }$s $S_{j,\text{ }k}$ are non increasing functions of
independent rv's $E_{j}.$ So they are associated. We then use the Lemma 3 of
Newman \cite{newman80} stated in Lemma \ref{lem1} to get%
\begin{equation*}
\left\vert Cov\left( \exp \left( -\gamma \sum_{h=i}^{k-1}E_{h}/h\right) ,%
\text{ }\exp \left( -\gamma \sum_{h=j}^{k-1}E_{h}/h\right) \right)
\right\vert
\end{equation*}%
\begin{equation*}
\leq Cov\left( \gamma \sum_{h=i}^{k-1}E_{h}/h,\text{ }\gamma
\sum_{h=j}^{k-1}E_{h}/h\right) ,
\end{equation*}%
where we use the one-value bound of $exp(-x)$. For $i\leq j,$%
\begin{equation}
Cov\left( \gamma \sum_{h=i}^{k-1}E_{h}/h,\text{ }\gamma
\sum_{h=j}^{k-1}E_{h}/h\right) =\mathbb{V}ar\left( \gamma
\sum_{h=j}^{k-1}E_{h}/h\right) =\gamma ^{2}\sum_{h=j}^{k-1}h^{-2}\leq \frac{%
\gamma ^{2}}{j},  \label{gc03-1}
\end{equation}%
the latter inequality is directly obtained by comparing $%
\sum_{h=j}^{k-1}h^{-2}$ and $\int_{j}^{k}x^{-2}dx.$ We get%
\begin{equation*}
|B_{k,2}|\leq \frac{1}{k^{1+\nu }}\sum_{L\leq i\neq j\leq k}\alpha
^{2}(k)\Delta f(j)\Delta f(i)Cov\left( \gamma \sum_{h=i}^{k-1}E_{h}/h,\text{ 
}\gamma \sum_{h=j}^{k-1}E_{h}/h\right)
\end{equation*}%
\begin{equation*}
\leq \frac{2\gamma ^{2}}{k^{1+\nu }}\sum_{L\leq i<j\leq k}\alpha
^{2}(k)\Delta f(j)\Delta f(i)/j
\end{equation*}%
\begin{equation}
=\frac{2\gamma ^{2}}{k^{1+\nu }}\alpha ^{2}(k)\sum_{j=L+1}^{k-1}\left[
\sum_{i=L}^{j-1}\Delta f(i)\right] \Delta f(j)\frac{1}{j}.  \label{gc03}
\end{equation}%
Finally, by using the techniques of (\ref{gc03-1}) and (\ref{gc03}), we get%
\begin{equation*}
C_{k}=\sum_{1\leq i\leq L-1}\sum_{L\leq j\leq k-1}\alpha ^{2}(k)\Delta
f(i)\Delta f(j)Cov(S_{i,\text{ }k},\text{ }S_{j,\text{ }k})
\end{equation*}%
\begin{equation}
\leq \frac{\alpha ^{2}(k)\gamma ^{2}}{k^{1+\nu }}\sum_{L\leq j\leq k-1}\left[
\sum_{1\leq i\leq L-1}\Delta f(i)\right] \Delta f(j)/j,  \label{gc04}
\end{equation}%
where $\left[ \sum_{1\leq i\leq L-1}\Delta f(i)\right] $ is a constant. By
putting together (\ref{gc01}), (\ref{gc02}), (\ref{gc03}) and (\ref{gc04}),
we get that assumptions (\ref{1}), (\ref{2}) and (\ref{2m}) entail (\ref%
{GCIP1}) in Theorem \ref{theo1}. We are going to check for (\ref{GCIP2})
now. We already noticed that the rv's $\alpha (k)\Delta f(q^{2}+i)(S_{k,%
\text{ }q^{2}+i}-s_{k,\text{ }q^{2}+i})$ are associated and partial sums of
associated rv's have non decreasing variances. Then for $j\leq 2q+1,$ we
have 
\begin{equation*}
\mathbb{V}ar\left( \sum_{i=1}^{j}\alpha (k)\Delta f(q^{2}+i)\left( S_{k,%
\text{ }q^{2}+i}-s_{k,\text{ }q^{2}+i}\right) \right)
\end{equation*}%
\begin{equation*}
\leq \mathbb{V}ar\left( \sum_{i=1}^{2q+1}\alpha (k)\Delta f(q^{2}+i)\left(
S_{k,\text{ }q^{2}+i}-s_{k,\text{ }q^{2}+i}\right) \right) .
\end{equation*}%
And (\ref{GCIP2}) becomes 
\begin{equation}
\sup_{k\geq 1}\sup_{(q^{2}+1)\leq k\leq (q+1)^{2}}\frac{1}{q^{(3-\delta )}}%
\mathbb{V}ar\left( \sum_{i=1}^{2q+1}\alpha (k)\Delta f(q^{2}+i)\left( S_{k,%
\text{ }q^{2}+i}-s_{k,\text{ }q^{2}+i}\right) \right) .  \label{ci20}
\end{equation}%
We fix $q$ but large enough to ensure $q^{2}\geq L,$ where $L$ is introduced
in (\ref{fadja01}). So (\ref{ci20}) is bounded by%
\begin{equation}
\sup_{(q^{2}+1)\leq k\leq (q+1)^{2}}\frac{1}{q^{(3-\delta )}}\mathbb{V}%
ar\left( \sum_{i=1}^{2q+1}\alpha (k)\Delta f(q^{2}+i)\left( S_{k,\text{ }%
q^{2}+i}-s_{k,\text{ }q^{2}+i}\right) \right) .  \label{ci21}
\end{equation}%
Now, we only have to show that 
\begin{equation*}
D=\sup_{(q^{2}+1)\leq k\leq (q+1)^{2}}\frac{1}{q^{(3-\delta )}}\mathbb{V}%
ar\left( \sum_{i=1}^{2q+1}\alpha (k)\Delta f(q^{2}+i)\left( S_{k,\text{ }%
q^{2}+i}-s_{k,\text{ }q^{2}+i}\right) \right)
\end{equation*}%
is bounded for $q^{2}\geq L.$ Let us split term in the brackets into%
\begin{equation*}
D=\frac{1}{q^{(3-\delta )}}\sum_{i=1}^{2q+1}\alpha ^{2}(k)\Delta
^{2}f(q^{2}+i)\mathbb{V}ar\left( S_{k,\text{ }q^{2}+i}\right)
\end{equation*}%
\begin{equation*}
+\frac{1}{q^{(3-\delta )}}\sum_{1\leq i\neq j\leq 2q+1}\alpha ^{2}(k)\Delta
f(q^{2}+i)\Delta f(q^{2}+j)Cov\left( S_{k,\text{ }q^{2}+i},\text{ }S_{k,%
\text{ }q^{2}+j}\right)
\end{equation*}%
\begin{equation*}
=:D_{1}+D_{2}.
\end{equation*}%
We have, by (\ref{fadja01}), 
\begin{equation*}
D_{1}=\frac{1}{q^{(3-\delta )}}\sum_{i=1}^{2q+1}\alpha ^{2}(k)\Delta
^{2}f(q^{2}+i)\mathbb{V}ar\left( S_{k,\text{ }q^{2}+i}\right)
\end{equation*}%
\begin{equation}
\leq (1/2)\frac{1}{q^{(3-\delta )}}\sum_{i=1}^{2q+1}\alpha ^{2}(k)\Delta
^{2}f(q^{2}+i)\left( \frac{q^{2}+i}{k}\right) ^{2\gamma }.  \label{ci22}
\end{equation}%
Now, we handle $D_{2}.$ We use again the techniques that lead to (\ref{gc03}%
) based on the Newman's Lemma to get, for $i\leq j,$ 
\begin{equation*}
\left\vert D_{2}\right\vert \leq \frac{1}{q^{(3-\delta )}}\sum_{1\leq i\neq
j\leq 2q+1}\alpha ^{2}(k)\Delta f(q^{2}+i)\Delta f(q^{2}+j)\mathbb{V}%
ar\left( \gamma \sum_{h=q^{2}+j}^{2q+1}E_{h}/h\right) .
\end{equation*}%
We remind, as in (\ref{gc03-1}), that 
\begin{equation*}
\mathbb{V}ar\left( \gamma \sum_{h=q^{2}+j}^{2q+1}E_{h}/h\right) \leq \gamma
^{2}/(q^{2}+j)
\end{equation*}%
and then%
\begin{equation*}
\left\vert D_{2}\right\vert \leq \frac{2\gamma ^{2}}{q^{(3-\delta )}}\alpha
^{2}(k)\sum_{1\leq i<j\leq 2q+1}\Delta f(q^{2}+i)\Delta f(q^{2}+j)\frac{1}{%
q^{2}+j}
\end{equation*}%
\begin{equation}
=\frac{2\gamma ^{2}}{q^{(3-\delta )}}\alpha ^{2}(k)\sum_{j=2}^{2q+1}\left[
\sum_{i=1}^{j-1}\Delta f(q^{2}+i)\right] \Delta f(q^{2}+j)\frac{1}{q^{2}+j}.
\label{ci23}
\end{equation}%
By putting together (\ref{ci22}) and (\ref{ci23}), we get that assumptions (%
\ref{3}) and (\ref{4}) entail (\ref{GCIP2}) in Theorem \ref{theo1}. We may
conclude that the strong law of large numbers holds for $S_{k}^{\ast }.$

\subsubsection{Special case for $f(j)=j^{\protect\tau }$}

We are going to check the conditions (\ref{1}), (\ref{2}), (\ref{2m}), (\ref%
{3}) and (\ref{4}) for the special function $f(j)=j^{\tau },$ $\tau >0.$ We
fix $L$ as indicated, consider $q\geq L$ and work with $k\geq q^{2}+1.$ We
notice that $\Delta f(j)$ is equivalent to $\tau j^{\tau -1}$ and $\Delta
f(q^{2}+j)$ is uniformly equivalent to $\tau j^{\tau -1}$ uniformly in $%
j\geq L.$ Here $\alpha (k)=k^{-(\tau -1)}.$ Then (\ref{1}) holds when%
\begin{equation*}
\sup_{k\geq L}\frac{\tau ^{2}}{k^{2\gamma +2\tau -1+\nu }}%
\sum_{j=L}^{k-1}j^{2\gamma +2\tau -2}
\end{equation*}%
is bounded. But if $2\gamma +2\tau -1=0,$ we get%
\begin{equation*}
\frac{1}{k^{\nu }}\sum_{j=L}^{k-1}j^{-1}\sim k^{-\nu }\log k\rightarrow 0
\end{equation*}%
and for $2\gamma +2\tau -1\neq 0,$ we get%
\begin{equation*}
\frac{1}{k^{2\gamma +2\tau -1+\nu }}\sum_{j=L}^{k-1}j^{2\gamma +2\tau
-2}\sim k^{-\nu }(2\gamma +2\tau -1)^{-1}
\end{equation*}%
and (\ref{1}) holds. (\ref{2}) holds with boundedness of 
\begin{equation*}
\sup_{k\geq L}\frac{1}{k^{2\tau -1+\nu }}\sum_{j=L+1}^{k-1}j^{2\tau -2}
\end{equation*}%
which is, for $2\tau -1\neq 0$ 
\begin{equation*}
\frac{1}{2\tau -1}k^{-\nu }\rightarrow 0,
\end{equation*}%
and for $2\tau =1$%
\begin{equation*}
k^{-\nu }\ln k\rightarrow 0.
\end{equation*}%
Next (\ref{2m}) is equivalent to the boundedness of 
\begin{equation*}
\frac{1}{k^{2\tau -1+\nu }}\sum_{j=L}^{k-1}j^{\tau -2},
\end{equation*}%
which is equivalent to the boundedness of $k^{-(\tau +\nu )}\log k,$ for $%
\tau -1=0$ and to that of $k^{-(\tau +\nu )}$ for $\tau -1\neq 0.$ Let us
now handle (\ref{3}) which is equivalent to the boundedness of 
\begin{equation*}
\frac{1}{q^{(3-\delta )}}\alpha ^{2}(k)\sum_{j=2}^{2q+1}\left[
\sum_{i=1}^{j-1}\Delta f(q^{2}+i)\right] \Delta f(q^{2}+j)\frac{1}{q^{2}+j}
\end{equation*}%
\begin{equation*}
\leq \frac{1}{q^{(3-\delta )}}\alpha ^{2}(k)\sum_{j=1}^{2q+1}\left[
\sum_{i=1}^{j-1}\Delta f(q^{2}+i)\right] \Delta f(q^{2}+j)\frac{1}{q^{2}+j},
\end{equation*}%
for $q$ large enough. We have to establish that%
\begin{equation*}
\sup_{k\geq 1}\sup_{(q^{2}+1)\leq k\leq (q+1)^{2}}\frac{1}{q^{(3-\delta )}}%
\frac{1}{k^{2\gamma +2\tau -2}}\sum_{j=1}^{2q+1}(q^{2}+j)^{2\gamma +2\tau
-2}<+\infty .
\end{equation*}%
If $2\gamma +2\tau -1\neq 0,$ then $\frac{1}{q^{(3-\delta )}}\frac{1}{%
k^{2\gamma +2\tau -2}}\sum_{j=1}^{k-(2k+1-q^{2})}(q^{2}+j)^{2\gamma +2\tau
-2}$ is bounded whenever 
\begin{equation*}
\frac{1}{q^{(3-\delta )}}\frac{1}{k^{2\gamma +2\tau -2}}\frac{k^{2\gamma
+2\tau -1}}{2\gamma +2\tau -1}=\frac{1}{2\gamma +2\tau -1}%
(k/q^{2})q^{-(1-\delta )}
\end{equation*}%
is bounded. And if $2\gamma +2\tau -1=0,$ $%
\sum_{j=1}^{k-(2k+1-q^{2})}(q^{2}+j)^{2\gamma +2\tau -2}$ is bounded along
with%
\begin{equation*}
\frac{k}{q^{(3-\delta )}}\log k\leq \left( k/q^{2}\right) q^{-(1-\delta
)}\log k.
\end{equation*}%
In both cases, $\left( k/q^{2}\right) q^{-(1-\delta )}\sim q^{-(1-\delta
)}\rightarrow 0$ as $k$ (and $q)$ goes to infinity. The proof is now
complete.\\

\noindent \textbf{ACKOWLEDGEMENT}
The first author thanks the \textbf{Programme de formation des formateurs} of USSTB who financed his stays in the LERSTAD of UGB while preparing his PhD dissertation. The authors aknowledge support from the \textbf{Réseau EDP - Modélisation et Contrôle}, of Western African Universities, that financed travel and accomodation of the second author while visiting USTTB in preparation of this work.

\end{document}